\documentclass[12pt]{article}
\usepackage[T2A]{fontenc}
\usepackage[cp1251]{inputenc}
\usepackage[english]{babel}
\usepackage{amsfonts}
\usepackage{amsmath}
\usepackage{amssymb}
\usepackage{graphicx}
\usepackage{graphics}
\usepackage[mathscr]{eucal}
\usepackage{calrsfs}

\newcommand{\blem}{\begin{lemma}$\!\!\!${\bf  .}~}
\newcommand{\elem}{\end{lemma}}
\newcommand{\bthm}{\begin{theorem}$\!\!\!${\bf  .}~}
\newcommand{\ethm}{\end{theorem}}
\newtheorem{note}{\hspace*{6mm}Remark}
\newcommand{\bnot}{\begin{note}$\!\!\!${\bf  .}~}
\newcommand{\enot}{\end{note}}

\newcommand{\sq}{\hfill$\square $}
\newcommand{\prooff}{\par{\it Proof. }\/}
\newcommand{\beqn}{\begin{eqnarray}}
\newcommand{\eeqn}{\end{eqnarray}}

\newtheorem{lemma}{$\phantom{.....}$Lemma}
\newtheorem{theorem}{$\phantom{.....}$Theorem}

\newcommand{\ris}{\refstepcounter{ris}\arabic{ris}}

\newcounter{ia}
\begin{document}

\thispagestyle{empty}
\newcounter{ris}
\setcounter{ris}{0} \vspace*{1mm}
\begin{center}\Large\bf
On an extreme two-point distribution\end{center}

\begin{center}\large Chebotarev V.I., Kondrik A.S., Mikhaylov K.V.\footnote{Computing Center of the
Far-Eastern Branch of the Russian Academy of Sciences
} \end{center}
\vskip3mm

Let  $\mathbb V$ be the class of all distribution function $F$ of
random variables with zero mean and unit variance. Denote
$$\varphi(t)=\frac{1}{\sqrt{2\pi}}\,e^{-t^2/2},\;\;\Phi(x)=\int_{-\infty}^x\varphi(t)\,dt,\;\;
\Delta(F)=\sup_{x\in\mathbb R}|F(x)-\Phi(x)|.$$

There was obtained in \cite[Lemma 12.3, p. 115]{Bh} (see also
\cite[Problem 3.269, p.~70]{Prohorov_U_U}) that
$
\sup\limits_{F\in{\mathbb
V}}\Delta(F)\le0.5416$.

We show in the present notice that in fact the proof in~\cite{Bh}
leads to a  bound, which is less than 0.541, and, moreover, this new
bound is reached.

To formulate the main statement we need the function
$$\Psi(x)=\frac{1}{1+x^2}-\Phi\left(-|x|\right).$$

{\bf Lemma}.
{\it The function $\Psi(x)$ is positive for all $x\in{\mathbb R}$,
and attains its maximum $C_\Phi=0.5409365\ldots\;$ in two points:
$x_\Phi=0.213105\ldots\;$ and $-x_\Phi$. In addition, $x_\Phi$ is
the root of the equation: $\;xe^{x^2/2}
(1+x^2)^{-2}=(8\pi)^{-1/2}$.}\vskip3mm

\prooff Let $x\ge0$. We have
\beqn\label{Psi'}\Psi^\prime(x)=\frac{-2x}{(1+x^2)^2}+\frac{1}{\sqrt{2\pi}}\;e^{-\frac{x^2}{2}}
 =-2e^{-\frac{x^2}{2}}\left(u(x)-\frac{1}{\sqrt{8\pi}}\right),\eeqn
 where $u(x)=\frac{xe^{x^2/2}}
{(1+x^2)^2}$.
 Since
$$(\ln u(x))^\prime=\frac{1}{x}+x-\frac{4x}{1+x^2}=\frac{(1-x^2)^2}{x(1+x^2)}
\quad\begin{cases}\le0,&\text{if}\quad x<0,\\
\ge0,&\text{if}\quad x>0,\end{cases}$$ then $u(x)$ increases, when
$x>0$ (there are similar computations in \cite{BenKir89} too).
With the help of computer we find that the equation
$u(x)-\frac{1}{\sqrt{8\pi}}=0$  has the root
$x_\Phi=0.21310518\ldots\;$, when  $x>0$, and, moreover,
$$C_\Phi:=\Psi(x_\Phi)=0.5409365\ldots\;.$$  Then it follows from
\eqref{Psi'} that $\Psi^\prime\ge0$ for $0< x\le x_\Phi$, and
$\Psi^\prime<0$ for $x> x_\Phi$. Hence,
$\max\limits_{x>0}\Psi(x)=C_\Phi$. Since  $\Psi$ is an even
function, the lemma is proved. For obviousness,
see~Fig.~\ref{sup-ris1} and Fig.~\ref{sup-ris2}.

\vskip6mm

\input epsf.tex
\hspace{-0.4cm} \epsfxsize=7cm \epsfbox{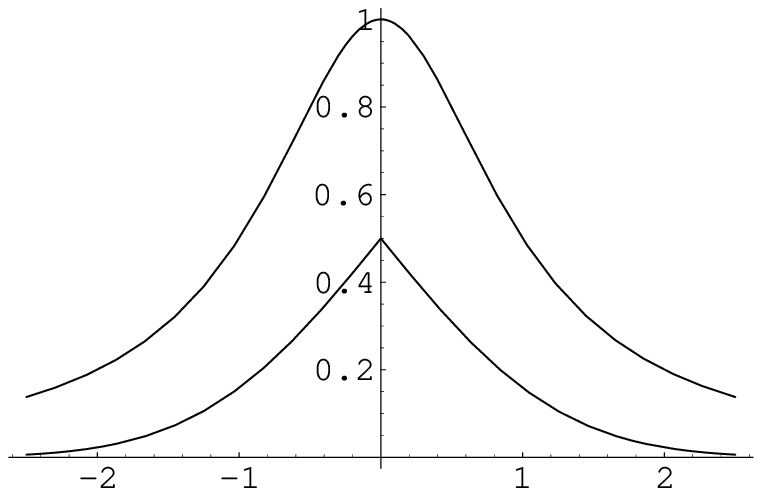} \epsfxsize=7cm
\quad\quad \epsfbox{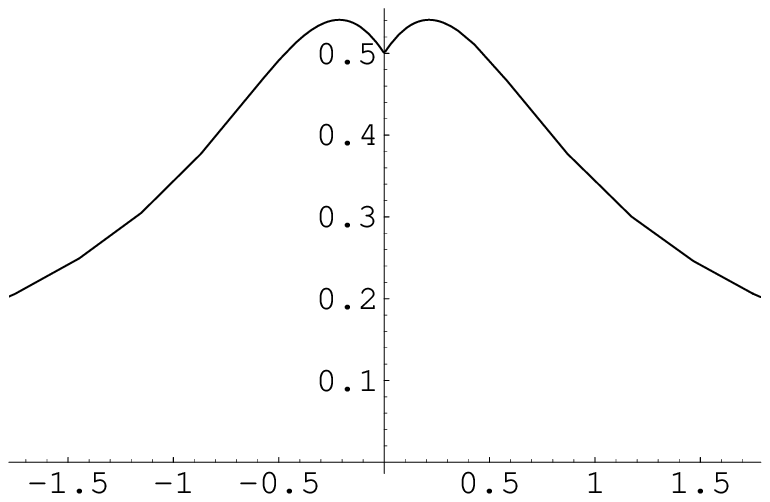} \vskip 0.5cm
\begin{tabular}{p{1cm}p{5cm}p{1mm}p{1cm}p{6cm}}
\mbox{\small\bf Fig. \ris\label{sup-ris1}.}&\hspace{2mm}\small The
functions $\frac{1}{1+x^2}$ \par\noindent \hspace{2mm} and
$\Phi(-|x|)$&& \mbox{\small\bf Fig.
\ris\label{sup-ris2}.}&\hspace{2mm} \mbox{\small The bihavior of the
function $\Psi(x)$}
\end{tabular}

 \sq

\bnot {\rm     It should
be noted that in \cite[Lemma 12.3, p. 115]{Bh} the less precise
values for $x_\Phi$ and $C_\Phi$ (in our notations) are indicated:
$x_\Phi=0.2135$, $C_\Phi=0.5416$. \sq}\enot

Define r.v. $X_{\Phi}$, taking two values,
$x_1=-\frac{1}{x_\Phi}=-4.692518\ldots\;$ and
$x_2=x_\Phi=0.213105\ldots\;$, with the probabilities
$p_1=\frac{x_\Phi^2}{1+x_\Phi^2}=0.043441\ldots\;$ and
$p_2=\frac{1}{1+x_\Phi^2}=0.956559\ldots\;$  respectively. It is not
hard to see  that the distribution functions of $X_\Phi$
and~$-X_\Phi$, denote them by $F_\Phi$ and~$\widetilde F_\Phi$
respectively, belong to the class $\mathbb V$.\vskip3mm

{\bf Theorem}.  {\it For every $F\in\mathbb V$,
\begin{equation}\label{main_ineq}
\Delta(F)\leq C_\Phi.
\end{equation}
The inequality \eqref{main_ineq} becomes the equality for
$F=F_\Phi$ and $F=\widetilde F_\Phi$.}

\bnot  {\rm The  following variant of the Berry -- Esseen
inequality, when the number of summands is equal to one, is proved
in \cite{BenKir89}: for every $\beta\ge1$
\begin{equation}\label{Ben-oc-1}\sup\limits_{F\in{\mathbb V}_\beta}\Delta(F)\le C_1\beta,\end{equation}
where $C_1=0.37035\ldots\;$, ${\mathbb V}_\beta$ is the subclass of
all such distribution functions from~${\mathbb V}$, that the third
absolute moment is equal to $\beta$. The bound \eqref{Ben-oc-1} is
unimprovable, provided that the right-hand side is the product of an
absolute constant and $\beta$.

Notice that if $\beta>C_\Phi/C_1\approx1.46$, the bound  \eqref{main_ineq} is more
precise with respect to~\eqref{Ben-oc-1}. But in the case
$1\le\beta<C_\Phi/C_1$ the inequality \eqref{Ben-oc-1} gives more
precise bounds for $\sup\limits_{F\in{\mathbb V}_\beta}\Delta(F)$,
than Theorem.\sq}\enot

{\it Proof of Theorem.} It is well-known, that if the distribution function $F$ of
r.v. $X$ belongs to $\mathbb V$, then \beqn\label{sup-nerav}{\bf
P}(X\le-x)\le\frac{1}{1+x^2}\quad\text{и}\quad {\bf P}(X\ge
x)\le\frac{1}{1+x^2},\quad x>0.\eeqn See, for instance, \cite[Lemma
12.3, p. 115]{Bh} and \cite[Problem 3.237, pp.~67,
239]{Prohorov_U_U}.

  Let $x
> 0$ be fixed. If $F(x)>\Phi(x)$, then, evidently,
$F(x)-\Phi(x)<0.5$. Consider the case $F(x)\le   \Phi(x)$. It is
easily seen that
$$
|F(x)-\Phi(x)|=\Phi(x)-F(x)=\Phi(x)-1+1-F(x)=-\int\limits_x^\infty
\varphi(t)\,dt+ {\bf P}(X \geq x) .$$ Using \eqref{sup-nerav}, we
obtain the inequality $ |F(x)-\Phi(x)|
 \leq \frac{1}{1+x^2} - \int\limits_x^\infty \varphi(t) \,dt=\Psi(x).$
It now follows from Lemma that
\begin{equation}\label{F-Phi<C} |F(x)-\Phi(x)|
 \le C_\Phi.\end{equation}
The inequality \eqref{F-Phi<C} for $x<0$ is deduced  similarly but
simpler. The bound \eqref{main_ineq} is proved.

Let now $F=F_\Phi$. Then we have

\begin{multline*}
\Phi(x_2)-F_\Phi(x_2)=\Phi(x_2) - p_1 = \Phi(x_2) - (1-p_2) \\= p_2
- (1 - \Phi(x_2)) = \frac{1}{1+x_2^2} - \int\limits_{x_2}^\infty
\varphi(t) \,dx=\Psi(x_2)=\Psi(x_\Phi)=C_\Phi.
\end{multline*}
For the sake of obviousness see Fig. \ref{F_phi}. Thus, the equality is attained in \eqref{main_ineq},
when
 $F=F_\Phi$.
\vspace*{-2mm}

 \begin{picture}(400,200)(-50,-60)
 {\small
 \put(0,30){\vector(1,0){270}}
 \put(210,29){$\mbox{\Large\bf.}$}
 \put(240,29){$\mbox{\Large\bf.}$}
 \put(210,20){$1$}
 \put(240,20){$2$}
 \put(120,29){$\mbox{\Large\bf.}$}
 \put(150,29){$\mbox{\Large\bf.}$}
 \put(112,20){$-2$}
 \put(142,20){$-1$}
 \put(60,29){$\mbox{\Large\bf.}$}
 \put(90,29){$\mbox{\Large\bf.}$}
 \put(30,29){$\mbox{\Large\bf.}$}
 \put(25,20){-5}
 \put(40,29){$\mbox{\Large\bf.}$}
 \put(40,20){$x_1$}
 \put(188,29){$\mbox{\Large\bf.}$}
 \put(188,20){$x_2$}
 \put(180,10){\vector(0,1){110}}
 \put(178,65){$\mbox{\Large\bf.}$}
 \put(178,100){$\mbox{\Large\bf.}$}
 \put(163,65){$0.5$}
 \put(172,100){$1$}

 \put(42.1,35){\line(1,0){148}}
 \put(190,101.5){\line(1,0){67}}
 \qbezier(180,65)(187,74)(195,80.5)
 \qbezier(180,65)(173,56)(165,49.5)
 \qbezier(195,80.5)(219,100)(243,100)
 \qbezier(165,49.5)(141,32)(117,32.5)

 \multiput(190,35)(0,5){14}{\line(0,7){1}}
 \put(160,53){$\left.\phantom{\displaystyle\int_{A^{B^B}}^{A^{B^B}}}\right\}C_\Phi$}
 \put(10,-10){{\small\bf Fig. \ris.\label{F_phi}} \small The distribution functions $\Phi(x)$
 and $F_\Phi$
} } \multiput(42.1,30)(0,2){3}{\line(0,7){1}}
\end{picture}\vspace*{-1cm}

\par\noindent   Symmetric considerations lead
 to the conclusion that for $F=\widetilde
F_\Phi$ the equality in \eqref{main_ineq} is attained   as well.\sq

\bnot{\rm The proof of \eqref{main_ineq}  differs from the proof of
Lemma  12.3~\cite{Bh} in fact only  by the deduction of the
extreme two-point distribution~$F_\Phi$, showing that the
bound~\eqref{main_ineq} is  unimprovable.\sq} \enot

{\bf Acknowledgments}. The work is partially supported
by the  grant RFBR 07-01-00054, and grant of
the Far-Eastern Branch of the Russian Academy of Sciences  06-III-A-01-003.

\label{lit}

\end{document}